\documentclass[10pt]{article}
\usepackage{amsmath, amsthm, amssymb }
\usepackage{float,epsfig}
\usepackage{color}
\usepackage{multirow}
\usepackage{cleveref}
\usepackage{braket}
\usepackage{subfigure}
\usepackage[utf8]{inputenc}
\usepackage{bbm}

\begin{document}

\newtheorem{theorem}{\bf Theorem}[section]
\newtheorem{proposition}[theorem]{\bf Proposition}
\newtheorem{definition}[theorem]{\bf Definition}
\newtheorem{corollary}[theorem]{\bf Corollary}
\newtheorem{example}[theorem]{\bf Example}
\newtheorem{remark}[theorem]{\bf Remark}
\newtheorem{lemma}[theorem]{\bf Lemma}
\newcommand{\nrm}[1]{|\!|\!| {#1} |\!|\!|}

\newcommand{\ba}{\begin{array}}
\newcommand{\ea}{\end{array}}
\newcommand{\dm}[1]{ {\displaystyle{#1} } }
\newcommand{\be}{\begin{equation}}
\newcommand{\ee}{\end{equation}}
\newcommand{\beano}{\begin{eqnarray*}}
\newcommand{\eeano}{\end{eqnarray*}}
\newcommand{\inp}[2]{\langle {#1} ,\,{#2} \rangle}
\def\bmatrix#1{\left[ \begin{matrix} #1 \end{matrix} \right]}
\def \noin{\noindent}
\newcommand{\evenindex}{\Pi_e}

\newcommand{\tb}[1]{\textcolor{blue}{ #1}}
\newcommand{\tm}[1]{\textcolor{magenta}{ #1}}
\newcommand{\tre}[1]{\textcolor{red}{ #1}}



\def \B{{\mathcal B}}
\newcommand{\C}{{\mathbb C}}
\def \H{{\mathcal H}}
\def \L{{\mathsf L}}
\def \O{{\mathcal O}}
\def \P{{\mathcal P}}
\def \Q{{\mathbb Q}}
\newcommand{\R}{{\mathbb R}}
\newcommand{\U}{{\mathrm U}}
\def \X{{\mathcal X}}
\def \Y{{\mathcal Y}}
\def \Z{{\mathcal Z}}

\def \PO{{\mathcal {PO}}}
\def \pf{{\bf Proof: }}
\def \lam{{\lambda}}

\def \sig{\sigma}



\title{A characterization of orthogonal permutative matrices of order 4}
\author{ Amrita Mandal\thanks{Department of Mathematics, IIT Kharagpur,Email: mandalamrita55@gmail.com}, \, and \,  Bibhas Adhikari\thanks{Corresponding author, Department of Mathematics, IIT Kharagpur,Email: bibhas@maths.iitkgp.ac.in}}
\date{}
\maketitle

\thispagestyle{empty}

\noindent{\bf Abstract.} Orthogonal matrices which are linear combinations of permutation matrices have attracted enormous attention in quantum information and computation. In this paper, we provide a complete parametric characterization of all complex, real and rational orthogonal permutative matrices of order $4.$ We show that any such matrix can always be expressed as a linear combination of up to four permutation matrices. Finally we determine several matrix spaces generated by linearly independent permutation matrices such that any orthogonal matrix in these spaces is always permutative or direct sum of orthogonal permutative matrices up to permutation of its rows and columns.\\




\section{Introduction}

Characterization of orthogonal matrices which can be expressed as linear combinations of permutation matrices is an unsolved problem in the literature. This problem was first considered by Kapoor in \cite{Kapoor1975} and he determined a necessary condition for a (real) linear combination of permutation matrices to be  an orthogonal matrix. Indeed, he proved that for a set $S$ of permutations of order $n,$ if a linear combination $\sum_{\sig\in S} \alpha_\sig P_\sig,$ $\alpha_\sig\in\R$ is an orthogonal matrix then $\sum_{\sig\in S} \alpha_\sig \in \{1,-1\}.$ Here $P_\sig$ denotes the permutation matrix associated with $\sig,$ that is, the $ij$ entry of $P_\sig$ is $1$ if $\sig(i)=j,$ otherwise it is $0.$ Later, Gibson proved that an orthogonal matrix over any field is a linear combination of permutation matrices if and only if each row and column sum of the matrix is $\pm 1,$ called generalized doubly stochastic (gds) matrix corresponding to $\pm 1$ \cite{gibson1980generalized}. However, explicit parametric representation of all gds matrices has remained to be determined.  


Parametric representation of orthogonal matrices which are also linear combinations of permutation matrices has recently attracted a lot of attention because they can have a proper quantum circuit representation \cite{klappenecker2003quantum}. Furthermore, orthogonal or unitary parametric matrices that are used as coin operators for discrete-time quantum walks play a crucial role in understanding the underlying quantum dynamics of the walks  \cite{watabe2008limit,vstefavnak2012continuous}. Recently, all orthogonal matrices of order $3$ over the fields of complex, real, and rational numbers that can be written as linear combinations of permutation matrices are characterized in \cite{sarkar2020}. One-parameter representations of such matrices are also provided. An interesting result about this characterization is that: an orthogonal matrix of order $3$ is a linear combination of permutation matrices if and only if it is a permutative matrix, that is, any row of such matrix is a permutation of any other row. We emphasize that the Grover diffusion matrix $G=\frac{2}{n}\boldsymbol{1}_n\boldsymbol{1}_n^T - I_n$, a standard coin operator for coined quantum walks is an orthogonal permutative matrix (OPM), where $\boldsymbol{1}_n$ is the all-one column vector of order $n$ and $I_n$ denotes the identity matrix of order $n.$ Thus the characterization of orthogonal permutative matrices of higher order has become of paramount interest since they can be used as coin operators for high dimensional quantum walks and hence a better understanding of quantum dynamics for generalized Grover diffusion matrices can be studied. In the forthcoming paper \cite{2dlocalization}, we investigate localization property of discrete-time quantum walks on two-dimensional lattices with coin operators as OPMs of order $4$ which are studied here.  

In this paper, we pursue the problem of algebraic characterization of the set of all OPMs of order $4$ and then investigate whether any orthogonal matrix of order $4$ which can be expressed as linear combination of permutation matrices belongs to this set. First, we derive symbolic representation of all OPMs of order $4$ over the field of complex numbers and we show that such a matrix can always be written as linear combination of up to four permutation matrices. However, contrary to the OPMs of order $3$, we establish that the set of OPMs of order $4$ does not form a group under matrix multiplication. Consequently, we determine chains of certain groups of OPMs of order $4$.  Then we develop a one-parameter representation of all real and rational OPMs of order $4$. 

Further, we produce an example of an orthogonal matrix of order $4$ which is a linear combination of permutation matrices but not permutative, and then we attempt to classify all such matrices by performing a search on the linear spaces of matrices generated by sets of linearly independent permutation matrices.  We prove that any linear combination of permutation matrices of order $4$ can be written as a linear combination of at most six permutative matrices, each of which is a linear combination of four pairwise Hadamard orthogonal permutation matrices, where two matrices are called Hadamard orthogonal if Hadamard product of them is the zero matrix. We also prove that there is no orthogonal matrix that can be expressed as a (non-trivial) linear combination of two distinct permutation matrices of order $4$. Next, we prove that any orthogonal matrix that can be expressed as a linear combination of three distinct permutation matrices, is always a direct sum of OPMs of orders $3$ and $1$ up to the permutations of its rows and its columns. Then, fixing a maximal set of linearly independent permutation matrices of order $4$, we determine direct sum of matrix spaces, which are generated by certain linearly independent permutation matrices, such that any orthogonal matrix $A$ in these spaces is always permutative or direct sum of OPMs up to the permutations of their rows and columns.

The paper is organised as follows. In Section 2, we provide a complete classification of complex, real and rational OPMs of order $4$ and show that any such matrix is always a linear combination of permutation matrices. Then we derive a one-parameter representation for the real and rational OPMs. In Section 3, we derive matrix spaces generated by permutation matrices such that orthogonal matrices in these spaces are always permutative or direct sum of OPMs up to permutation of rows and columns.  


\section{Orthogonal permutative matrices of order $4$}\label{Sec:2}

In this section we characterize the set of all OPMs of order $4$ over the field of complex, real and rational numbers. Recall that a matrix is called permutative if any of its row is a permutation of any other row  \cite{paparella2015realizing}. Thus, such a matrix of order $n$ has $n$ parameters in its symbolic form.  Without loss of generality, a permutative matrix of order $4$ can be written in the symbolic form
\begin{eqnarray} \label{permutative:A}
A(\mathbf{x}; P,Q,R) = 
\bmatrix{\mathbf{x}\\ \mathbf{x}P\\ \mathbf{x}Q\\ \mathbf{x}R}
\end{eqnarray}
where $\mathbf{x} =\bmatrix{x & y & z & w}$
is a symbolic row vector with $x, y, z, w\in \C$ and $P,Q,R \in \P_4$, the group of permutation matrices of order $4$ \cite{pdspaper}. 


Let $\O\P_4$ denote the set of all OPMs over the field of complex numbers of order $4$. Then obviously $\P_4 \subset \O\P_4$. We denote $$1\oplus \P_3 = \left\{ \bmatrix{1 & \textbf{0}^T \\ \textbf{0} & P}\in \P_4 : P\in \P_3\right\},$$ where $\textbf{0}=\bmatrix{0 & 0 & 0}^T$ and $\P_3$ denotes the group of all permutation matrices of order $3$. Then the following theorem characterizes all matrices in $\O\P_4.$

\begin{theorem} \label{thm1} 
A matrix $A \equiv A(\mathbf{x}; P,Q,R)$ given by equation (\ref{permutative:A}) is an OPM if and only if $A\in \X\cup \Y \cup \Z$ where
 \beano  \X &=& \left\{ \overline PM^{\pm}_{x,z},\overline PN^{\pm}_{z,x}:  x^2+z^2\mp z=0, x,z\in \C\right\} \\
   \Y &=& \left\{ \overline PP_{(23)}M^{\pm}_{x,z}P_{(23)},\overline PP_{(23)}N^{\pm}_{z,x}P_{(23)} : x^2+z^2\mp z=0, x, z\in\C\right\} \\
\Z &=& \left\{ \overline PP_{(24)}M^{\pm}_{x,z}P_{(24)},\overline PP_{(24)}N^{\pm}_{z,x}P_{(24)} : x^2+z^2\mp z=0, x, z\in\C\right\},
\,\mbox{with}
\eeano 
\beano &&  M^{\pm}_{x,z}= \bmatrix{A_x & B^{\pm}_z \\ B_z^{\pm} & -A_x}, N^{\pm}_{x,z}=\bmatrix{B^{\pm}_x & A_z \\ A_z & FB^{\pm}_x}, A_t=\bmatrix{t & -t \\ -t & t},  B^{\pm}_t= \bmatrix{t & \pm 1-t \\ \pm 1-t & t}  \\ &&  F=\bmatrix{0 & 1\\ 1 & 0},   t\in\{x,z\},\,\, \mbox{and} \,\, \overline P \in 1\oplus \P_3. \eeano
\end{theorem} 

\pf
The ‘if’ part is obvious and easy to check. To prove the ‘only if’ part consider the
following cases. First assume that the symbolic OPM $A$ has
no repetition of entries in any of the columns. Besides, since rows and columns are orthogonal, none of $ P,Q,R$ are equal to each other. Then $A$ can presume one of the following forms:
 \begin{equation}\label{eqn:matrix x,y,z}\mathsf{X} = \bmatrix{x&y&z&w\\
y&x&w&z\\
z&w&y&x\\
w&z&x&y}, \,\, \mathsf{Y} =  \bmatrix{x&y&z&w\\
y&z&w&x\\
z&w&x&y\\
w&x&y&z}, \,\, \mathsf{Z} = \bmatrix{x&y&z&w\\
y&w&x&z\\
z&x&w&y\\
w&z&y&x}\end{equation} for some $x, y, z, w\in\C.$ 

Then $\mathsf{X}^T\mathsf{X}=I_4$ provides the polynomial equations $x^2 + y^2 + z^2 + w^2 = 1,$ $x y + z w=0$ and $x z+y w+y z+z x=0,$ which further imply $(x+y+z+w)^2= 1.$ Thus  the quadruple $(x, y, z, w)$ must be zeros of the system of polynomial equations 
\begin{equation}\label{eqn:variety x}\begin{cases}
x + y + z + w=\pm1\\
x y + z w=0\\
x^2 + y^2 + z^2 + w^2 = 1
\end{cases} \Rightarrow \begin{cases}
x + y =0 \\
x^2 + z^2 \mp z = 0\\
z+w \mp 1 =0
\end{cases} \mbox{or} \begin{cases}
z + w =0\\
x^2 + z^2 \mp x=0\\
x+y \mp 1 =0.
\end{cases}\end{equation}
  Therefore, each of the four system of equations gives rise to the following set of matrices obeying the pattern $\mathsf{X}$: \begin{eqnarray}  \X_1 = \left\{ \bmatrix{A_x & B^+_z \\ B_z^+ & -A_x}:x^2+z^2-z=0\right\},
&& \X_2 = \left\{ \bmatrix{A_x & B^-_z \\ B^-_z & -A_x}:x^2+z^2+z =0\right\},\nonumber\\
 \X_3 =\left\{\bmatrix{B^+_x & A_z \\ A_z & FB^+_x} : x^2+z^2-x=0  \right\}, &&
\X_4 = \left\{\bmatrix{B_x^- & A_z \\ A_z & FB^-_x} : x^2+z^2+x =0 \right\}.\label{eqn:x1to4}\end{eqnarray}  Hence, \beano&\X_1=\{M^+_{x,z}:x^2+z^2-z =0\},\,\,\X_2=\{M^-_{x,z}:x^2+z^2+z =0\},\\
&\X_3=\{N^+_{x,z}:x^2+z^2-x =0\},\,\,\X_4=\{N^-_{x,z}:x^2+z^2+x =0\},\eeano where $M^{\pm}_{x,z},N^{\pm}_{x,z}, A_t, B_t, t\in\{x, z\}$ are defined in the statement of the theorem.

Similarly, the set of polynomial equations given by $\mathsf{Y}^T\mathsf{Y}=I_4$ are 
\begin{equation}\label{eqn:variety y}\begin{cases}
x + y + z + w = \pm 1\\
x z + y w=0\\
x^2 + y^2 + z^2 + w^2 = 1
\end{cases} \Rightarrow \begin{cases}
x + z =0 \\
x^2 + y^2 \mp y = 0\\
y+w \mp 1 =0
\end{cases} \mbox{or} \begin{cases}
y + w =0\\
x^2 + y^2 \mp x=0\\
x+z \mp 1 =0.
\end{cases}\end{equation}
  Thus each of the four system of equations gives rise to the following sets of matrices obeying the pattern $\mathsf{Y}$: \begin{eqnarray}  \Y_1 &=& \left\{ P_{(23)}\bmatrix{A_x & B^+_y \\ B_y^+ & -A_x}P_{(23)}:x^2+y^2-y=0\right\}, \nonumber\\
 \Y_2 &=& \left\{ P_{(23)}\bmatrix{A_x & B^-_y \\ B^-_y & -A_x}P_{(23)}:x^2+y^2+y =0\right\},\nonumber\\
 \Y_3 &=& \left\{P_{(23)}\bmatrix{B^+_x & A_y \\ A_y & FB^+_x}P_{(23)} : x^2+y^2-x=0  \right\}, \nonumber\\
\Y_4 &=& \left\{P_{(23)}\bmatrix{B_x^- & A_y \\ A_y & FB^-_x}P_{(23)} : x^2+y^2+x =0 \right\}.\label{eqn:y1to4}\end{eqnarray} 

Note that,
$A=P_{(23)}M^+_{x,y}P_{(23)}$ if $A\in\Y_1;$ $A=P_{(23)}M^-_{x,y}P_{(23)}$ if $A\in\Y_2;$ $A=P_{(23)}N^+_{x,y}P_{(23)}$ if $A\in\Y_3;$ and $A=P_{(23)}N^-_{x,y}P_{(23)}$ if $A\in\Y_4$ where $(x,y)$ satisfies the respective constraint as given in equation (\ref{eqn:y1to4}). 

Finally, the orthogonality condition $\mathsf{Z}^T\mathsf{Z}=I_4$ provides the system of polynomial equations
\begin{equation}\label{eqn:variety z}\begin{cases}
x + y + z + w = \pm 1\\
x w + y z=0\\
x^2 + y^2 + z^2 + w^2 = 1
\end{cases} \Rightarrow \begin{cases}
x + w =0 \\
x^2 + y^2 \mp y = 0\\
y+z \mp 1 =0
\end{cases} \mbox{or} \begin{cases}
y + z =0\\
x^2 + y^2 \mp x=0\\
x+w \mp 1 =0.
\end{cases}\end{equation}
  Thus each of the four system of equations gives rise to the following set of matrices obeying the pattern $\mathsf{Z}$: \begin{eqnarray}  \Z_1 = \left\{ P_{(24)}\bmatrix{A_x & B^+_y \\ B_y^+ & -A_x}P_{(24)}:x^2+y^2-y=0\right\},
& \Z_2 = \left\{ P_{(24)}\bmatrix{A_x & B^-_y \\ B^-_y & -A_x}P_{(24)}:x^2+y^2+y =0\right\},\nonumber\\
 \Z_3 =\left\{P_{(24)}\bmatrix{B^+_x & A_y \\ A_y & FB^+_x}P_{(24)} : x^2+y^2-x=0  \right\}, &
\Z_4 = \left\{P_{(24)}\bmatrix{B_x^- & A_y \\ A_y & FB^-_x}P_{(24)} : x^2+y^2+x =0 \right\}.\label{eqn:z1to4}\end{eqnarray} 
Besides, $A=P_{(24)}M^+_{x,y}P_{(24)}$ if $A\in \Z_1;$ $A=P_{(24)}M^-_{x,y}P_{(24)}$ if $A\in\Z_2;$ $A=P_{(24)}N^+_{x,y}P_{(24)}$ if $A\in\Z_3;$ and $A=P_{(24)}N^-_{x,y}P_{(24)}$ if $A\in\Z_4.$ 

Next, consider the symbolic OPMs in which one
entry is repeated in at least one column i.e. the case when $P,Q,R$ are not chosen from any of the collections $$\{P_{(12)(34)},P_{(1324)},P_{(1423)}\}, \{P_{(13)(24)},P_{(1234)},P_{(1432)}\}, \{P_{(14)(23)},P_{(1243)},P_{(1342)}\}.$$ Then it follows that any such matrix $A$ belongs to either of the sets $\X=\cup_{k=1}^4 \X_k, \Y=\cup_{k=1}^4 \Y_k, \Z=\cup_{k=1}^4 \Z_k$ for certain values of $x, y, z, w$. In particular, a straightforward calculation shows that any symbolic OPM whose entries do not follow the pattern of entries of $\mathsf{X}, \mathsf{Y}, \mathsf{Z}$ are  
$\pm P_{\tau}$ or $\pm(\frac{1}{2}J_4-P_{\tau})$ for some permutation $\tau $ and $J_4=\boldsymbol{1}_4\boldsymbol{1}_4^T.$ Then the desired result follows from the fact that permutation of the rows preserves the orthogonality and permutative property of a matrix.  $\hfill{\square}$

The following corollary provides the determinant of all OPMs of order $4.$
\begin{corollary}
Let $A\in \O\P_4.$ Then $\det(A)=1$ if $A\in \X_j \cup \Y_j\cup \Z_j, j=1, 2$, and $\det(A)=-1$ if $A\in \X_j\cup \Y_j\cup \Z_j, j=3, 4,$ where $\X_j, \Y_j, \Z_j, j=1, 2, 3, 4$ are given by equations (\ref{eqn:x1to4}), (\ref{eqn:y1to4}), (\ref{eqn:z1to4}) respectively.
\end{corollary}
\pf From (\ref{eqn:x1to4}) we obtain, $$\det(A)=\begin{cases} 4(x^2+z^2-z)+1 \, \mbox{if} \, A\in \X_1 \\ 4(x^2+z^2+z)+1 \, \mbox{if} \, A\in \X_2 \\ -4(x^2+z^2-x)-1  \, \mbox{if} \, A\in \X_3 \\ -4(x^2+z^2+x)-1 \, \mbox{if} \, A\in \X_4.
\end{cases}$$ Then employing the conditions on the variables $x,y,z$ which define the sets of matrices $\X_j, j=1,2,3,4$ the desired result follows. Similarly, the desired results for $\Y_j, \Z_j$ $j=1,2,3,4$ follow from equations (\ref{eqn:y1to4}) and (\ref{eqn:z1to4}). $\hfill{\square}$


Now in the following remark we provide a characterization of OPMs of order $4$ in terms of linear combinations of permutation matrices. The result follows from equation (\ref{eqn:matrix x,y,z}) and Theorem \ref{thm1}.
\begin{remark} \label{remark1} Any OPM $A$ of order $4$ can be written as linear combination of permutation matrices as follows:
\begin{equation}  \label{PO}
     \overline{P}^TA =\begin{cases}
     x P_{(34)} + y P_{(12)} + z P_{(13)(24)}+w P_{(14)(23)}, \,\, (x,y,z,w) \,\mbox{satisfies equation} \, (\ref{eqn:variety x}) \, \mbox{when} \,\, A\in\X \\
      x P_{(24)} + y P_{(12)(34)} + z P_{(13)} +w P_{(14)(23)}, \,\, (x,y,z,w) \,\mbox{satisfies equation} \, (\ref{eqn:variety y}) \, \mbox{when} \,\, A \in \Y \\
      x P_{(23)}+y P_{(12)(34)} +z P_{(13)(24)}+w P_{(14)}, \,\, (x,y,z,w) \,\mbox{satisfies equation} \, (\ref{eqn:variety z})\, \mbox{when} \,\, A \in \Z 
     \end{cases} 
\end{equation}
and $\overline{P}\in 1\oplus\P_3.$
\end{remark}

Let us emphasize that one of the motivations for the characterization of OPMs is to generalize Grover matrix which is used to define Grover quantum walks \cite{inui2004, vstefavnak2008recurrence, vstefavnak2010full}. The Grover matrix of order $4$ is given by  \begin{equation}\label{def:gm}G = \bmatrix{-\frac{1}{2}&\frac{1}{2}&\frac{1}{2}&\frac{1}{2}\\\frac{1}{2}&-\frac{1}{2}&\frac{1}{2}&\frac{1}{2}\\ \frac{1}{2}&\frac{1}{2}&-\frac{1}{2}&\frac{1}{2}\\ \frac{1}{2}&\frac{1}{2}&\frac{1}{2}&-\frac{1}{2}}.\end{equation} 

Then it can be seen that 
$$ G =\begin{cases} 
P_{(34)} \left(xP_{(34)} +yP_{(12)} + zP_{(13)(24)} + w P_{(14)(23)}\right) \in \X \\
P_{(24)}\left( x P_{(24)} + y P_{(12)(34)} + z P_{(13)} +w P_{(14)(23)}\right) \in \Y \\
P_{(23)} \left(x P_{(23)}+y P_{(12)(34)} +z P_{(13)(24)}+w P_{(14)} \right) \in \Z
\end{cases}
 $$ where $x=-\frac{1}{2}, y=\frac{1}{2}= z=w.$ In particular, it follows that \begin{equation}\label{eqn:groverxyz} G \in  P_{(34)}\X_1\cap P_{(24)}\Y_1\cap P_{(23)}\Z_1\end{equation} where \beano
    \X_1&=&\left\{\bmatrix{x&-x&z&1-z\\-x&x&1-z&z\\z&1-z&-x&x\\1-z&z&x&-x}:x^2+z^2-z=0\right\},\\
 \Y_1&=&\left\{ \bmatrix{x &y&-x&1-y\\y&-x&1-y&x\\-x&1-y&x&y\\1-y&x&y&-x}:x^2+y^2-y=0 \right\},\\
    \Z_1&=&\left\{\bmatrix{x&y&1-y&-x\\y&-x&x&1-y\\1-y&x&-x&y\\-x&1-y&y&x}:x^2+y^2-y=0\right\}.
 \eeano

Thus the real matrices in $P_{(34)}\X_1, P_{(24)}\Y_1, P_{(23)}\Z_1$ can be considered as the continuous deformations of the Grover matrix, and hence Grover walks can be generalized by considering the coin operators as matrices from these sets. The localization property of such quantum walks on two-dimensional lattices is analyzed in \cite{2dlocalization}.


Now we provide a list of chains of groups of OPMs of order $4$ in the following theorem. 

\begin{theorem}\label{prop:complex gp} 
The following are chains of groups of complex orthogonal matrices.
\begin{enumerate}
    \item  $\begin{array}{ll}
    & \{ I \} \leq {P_{(34)} \X_3} \leq {P_{(34)} \X_3 \cup P_{(34)} \X_j} \leq {P_{(34)} \X_3 \cup P_{(34)}\X_j \cup \X_3 \cup \X_j} \leq {\O_4},  \\  &  \{I\} \leq P_{(34)}\X_3 \leq P_{(34)}\X_3 \cup \X_j\leq  P_{(34)} \X_3 \cup P_{(34)}\X_j \cup \X_3 \cup \X_j \leq \O_4  
  \end{array}$
  \item  $\begin{array}{ll}
   &  \{ I \} \leq P_{(24)}\Y_3 \leq P_{(24)}\Y_3 \cup P_{(24)}\Y_j \leq P_{(24)}\Y_3 \cup P_{(24)}\Y_j \cup \Y_3 \cup \Y_j \leq \O_4,  \\ 
   & \{I\} \leq  P_{(24)}\Y_3 \leq P_{(24)}\Y_3 \cup \Y_j\leq  P_{(24)}\Y_3 \cup P_{(24)}\Y_j \cup \Y_3 \cup \Y_j \leq \O_4 
  \end{array}$
  \item  $\begin{array}{ll}
   & \{ I \} \leq  P_{(23)}\Z_3 \leq P_{(23)}\Z_3 \cup P_{(23)}\Z_j \leq P_{(23)}\Z_3 \cup P_{(23)}\Z_j \cup \Z_3 \cup \Z_j \leq \O_4,  \\ 
   & \{ I \} \leq  P_{(23)}\Z_3 \leq P_{(23)}\Z_3 \cup \Z_j\leq  P_{(23)}\Z_3 \cup P_{(23)}\Z_j \cup \Z_3 \cup \Z_j \leq \O_4 
  \end{array}$ 
\end{enumerate}
    where $j=1,2,3,4$ and $\O_4$ denotes the group of complex orthogonal matrices of order $4$.
 
\end{theorem}

\pf
First we prove that
$P_{(34)}\X_3 \cup P_{(34)}\X_j \cup \X_3 \cup \X_j$ are  complex orthogonal matrix groups for $j=1,2,3,4$. Clearly $I \in P_{(34)}\X_3 \cup P_{(34)}\X_j \cup \X_3 \cup \X_j$. 
If $A \in P_{(34)}\X_3$ then $A^T \in P_{(34)}\X_3 $  follows by exchanging the role of $z$ and $-z$. Similarly, $A^T \in P_{(34)}\X_j$ if $A \in P_{(34)} \X_j$ for $j=1,\ldots,4$. 
Since $\X_j$ and $\X_3$ contain complex symmetric matrices, obviously $A^T \in \X_j$ if $A \in \X_j$, and $A^T \in \X_3$ if $A \in \X_3$. Hence, $P_{(34)}\X_3 \cup P_{(34)}\X_j \cup \X_3 \cup \X_j,j=1,\ldots,4$ is closed under inverses.

Let $A = x_1 I + y_1 P_{(12)(34)} + z_1 P_{(1324)}+w_1 P_{(1423)},B = x_2 I + y_2 P_{(12)(34)} + z_2 P_{(1324)}+w_2 P_{(1423)} \in P_{(34)}\X_3 \cup P_{(34)}\X_j,j=1,\ldots,4,$ where $(x_i, y_i,z_i,w_i),i=1,2$ satisfy the equation given by (\ref{eqn:variety x}) accordingly. 
Then $AB = x_3 I + y_3 P_{(12)(34)} + z_3 P_{(1324)}+w_3 P_{(1423)},$ where
\begin{align*}
    &x_3=x_1x_2+y_1y_2+z_1w_2+w_1z_2,y_3=x_1y_2+y_1x_2+z_1z_2+w_1w_2,\\
    &z_3=x_1z_2+y_1w_2+z_1x_2+w_1y_2,w_3=x_1w_2+y_1z_2+z_1y_2+w_1x_2.
\end{align*}

Now note that $x_j+y_j+z_j+w_j=1,$  if $A\,\mbox{or}\,B\in P_{(34)}\X_1 \cup P_{(34)}\X_3$, and $x_j+y_j+z_j+w_j=-1$ if $A\,\mbox{or}\,B\in P_{(34)}\X_2 \cup P_{(34)}\X_4, j=1,2.$ 



Thus $x_3+y_3+z_3+w_3=(x_1+y_1+z_1+w_1)(x_2+y_2+z_2+w_2)$ yields  $x_3+y_3+z_3+w_3=1$ if $A,B \in P_{(34)}\X_1 \cup P_{(34)}\X_3~ \text{or}~ P_{(34)}\X_2 \cup P_{(34)}\X_4$, and $x_3+y_3+z_3+w_3=-1$ otherwise. Now $$x_3^2+y_3^2+z_3^2+w_3^2={\bmatrix{x_1\\y_1\\z_1\\w_1}}^T\bmatrix{a&b&c&c\\b&a&c&c\\c&c&a&b\\c&c&b&a}\bmatrix{x_1\\y_1\\z_1\\w_1},$$
where $a=x_2^2+y_2^2+z_2^2+w_2^2=1,b=2x_2y_2+2z_2w_2=0$ and $c=(x_2+y_2)(z_2+w_2)=0$ hold, since $B$ is orthogonal i.e. $(x_2,y_2,z_2,w_2)$ satisfies $x_2^2 + y_2^2 + z_2^2 + w_2^2 = 1, x_2 y_2 + z_2 w_2=0,x_2 z_2+y_2 w_2+y_2 z_2+w_2 x_2=0$. Consequently,  $x_3^2+y_3^2+z_3^2+w_3^2=x_1^2+y_1^2+z_1^2+w_1^2=1.$ Also,
\begin{align} \label{sum_x_3+y_3}
\begin{split}
  &x_3+y_3=(x_1+y_1)(x_2+y_2)+(z_1+w_1)(z_2+w_2) \\ &z_3+w_3=(x_1+y_1)(z_2+w_2)+(z_1+w_1)(x_2+y_2). 
  \end{split}
\end{align} 
At first, if $A,B\in P_{(34)}\X_1,$ then $x_i+y_i=0,z_i+w_i=1$ for $i=1,2.$ So that from (\ref{sum_x_3+y_3}) we have $x_3+y_3=1,z_3+w_3=0,$ where $x_3^2+y_3^2+z_3^2+w_3^2=1,$ and thus by (\ref{eqn:x1to4}) we get $AB\in P_{(34)}\X_3.$ If $A \in P_{(34)}\X_1$ and $B\in P_{(34)}\X_3,$ then clearly $x_1+y_1=z_2+w_2=0,x_2+y_2=z_1+w_1=1.$ Hence we get $x_3+y_3=0,z_3+w_3=1.$ Thus by (\ref{eqn:x1to4}) we have $AB\in P_{(34)}\X_3.$
Similarly, it can be done for other cases and we obtain
$AB \in P_{(34)}\X_3$ if $A,B\in P_{(34)}\X_j$ and $AB\in P_{(34)}\X_j$ if either $A$ or $B\in P_{(34)}\X_j,$ $j=1,\ldots,4.$
Thus finally $AB \in P_{(34)}\X_3 \cup P_{(34)}\X_j$ for $j=1,\ldots,4.$

Next, let $A = x_1 P_{(34)} + y_1 P_{(12)} + z_1 P_{(13)(24)}+w_1 P_{(14)(23)},B = x_2 P_{(34)} + y_2 P_{(12)} + z_2 P_{(13)(24)}+w_2 P_{(14)(23)} \in \X_3 \cup \X_j$ where $(x_i, y_i,z_i,w_i),i=1,2$ satisfies the equations in (\ref{eqn:variety x}) and $j=1,\ldots,4$. Then $AB = x_3 I+ y_3 P_{(12)(34)} + z_3 P_{(1324)}+w_3 P_{(1423)},$ where
\begin{align*}
    &x_3=x_1x_2+y_1y_2+z_1z_2+w_1w_2,y_3=x_1y_2+y_1x_2+z_1w_2+w_1z_2,\\
    &z_3=x_1z_2+y_1w_2+z_1y_2+w_1x_2,w_3=x_1w_2+y_1z_2+z_1x_2+w_1y_2.
\end{align*}

Note that $x_j+y_j+z_j+w_j=1$ if $A$ or $B\in \X_1 \cup \X_3,$ and $x_j+y_j+z_j+w_j=-1$ if $A$ or $B\in \X_2 \cup \X_4,j=1,2.$
Hence $x_3+y_3+z_3+w_3=1$ if $A,B \in \X_1 \cup \X_3~ \text{or}~ \X_2 \cup \X_4,$ and $x_3+y_3+z_3+w_3=-1$ otherwise, since $x_3+y_3+z_3+w_3=(x_1+y_1+z_1+w_1)(x_2+y_2+z_2+w_2).$  
Then as above, $x_3^2+y_3^2+z_3^2+w_3^2=1.$ 
Further, $x_3+y_3$ and $z_3+w_3$  have same expressions as that given in  (\ref{sum_x_3+y_3}). 

If $A,B\in \X_2,$ then $x_1+y_1=x_2+y_2=0,w_1+z_1=w_2+z_2=-1.$ So that $x_3+y_3=1$ and $z_3+w_3=0.$ Finally, $x_3^2+y_3^2+z_3^2+w_3^2=1$ implies $x_3^2+z_3^2-x_3=0.$ Thus $AB\in P_{(34)}\X_3.$ Now if $A\in \X_2$ and $B\in \X_3,$ then $x_1+y_1=z_2+w_2=0,w_1+z_1=x_2+y_2=1,$ which yield $x_3+y_3=0,z_3+w_3=-1.$ Now $x_3^2+y_3^2+z_3^2+w_3^2=1$ implies $x_3^2+z_3^2+z_3=0.$ Thus by (\ref{eqn:x1to4}) we have $AB\in P_{(34)}\X_2.$ The other cases can be done similarly and we obtain
$AB \in P_{(34)}\X_3~$ if $A,B\in \X_j$ and $AB\in P_{(34)}\X_j$ if either $A$ or $B\in \X_j$, $j=1,\ldots,4.$
Thus finally $AB \in P_{(34)}\X_3 \cup P_{(34)}\X_j$ for $j=1,\ldots,4.$

Now suppose $A = x_1 I + y_1 P_{(12)(34)} + z_1 P_{(1324)}+w_1 P_{(1423)}\in P_{(34)}\X_j$ and $B = x_2 P_{(34)} + y_2 P_{(12)} + z_2 P_{(13)(24)}+w_2 P_{(14)(23)} \in \X_k,$ $j,k\in\{1,\ldots,4\}$ and $(x_i,y_i,z_i,w_i),i=1,2$ are given by equation (\ref{eqn:variety x}) accordingly. Then $AB=x_3P_{(34)}+y_3P_{(12)}+z_3P_{(13)(24)}+w_3P_{(14)(23)},$
where 
\begin{align*}
    &x_3=x_1x_2+y_1y_2+z_1z_2+w_1w_2,y_3=x_1y_2+y_1x_2+z_1w_2+w_1z_2,\\
    &z_3=x_1z_2+y_1w_2+z_1y_2+w_1x_2,w_3=x_1w_2+y_1z_2+z_1x_2+w_1y_2.
\end{align*}
Thus by similar arguments as the above two cases, we obtain: 
$x_3+y_3+z_3+w_3=1$ if $A\in P_{(34)} \X_1 \cup P_{(34)}\X_3, B\in \X_1 \cup \X_3$ or $A\in P_{(34)} \X_2 \cup P_{(34)}\X_4, B\in \X_2 \cup \X_4$, and $x_3+y_3+z_3+w_3=-1$ otherwise. Then it can be checked that  $x_3^2+y_3^2+z_3^2+w_3^2=1$ and expressions of $x_3+y_3$ and $z_3+w_3$ are given by (\ref{sum_x_3+y_3}). Further, if $A\in P_{(34)}\X_1$ and $B\in \X_2,$ then $x_1+y_1=x_2+y_2=0,w_1+z_1=-(w_2+z_2)=1.$ So that $x_3+y_3=-1$ and $z_3+w_3=0.$ Finally, $x_3^2+y_3^2+z_3^2+w_3^2=1$ yields $x_3^2+z_3^2+x_3=0.$ Hence by (\ref{eqn:x1to4}) we write $AB\in \X_4.$ Similarly, the other cases follow and we get $AB \in \X_3~\text{if}~A\in P_{(34)}\X_j~\text{and}~B\in \X_j; AB\in \X_4~\text{if}~A \in P_{(34)}\X_j~\text{and}~B\in \X_k,j\neq k,\{j,k\}\in\{\{1,2\},\{2,1\},\{3,4\},\{4,3\}\}; AB \in \X_1~\text{if}~A \in P_{(34)}\X_j~\text{and}~B\in \X_k,j\neq k,\{j,k\}\in\{\{1,3\},\{3,1\},\{2,4\},\{4,2\}\};$ and $AB\in \X_2~\text{if}~A \in P_{(34)}\X_j~\text{and}~B\in \X_k,j\neq k,\{j,k\}\in\{\{1,4\},\{4,1\},\{2,3\},\{3,2\}\}.$


Thus considering all the above cases we conclude that $P_{(34)}\X_3,P_{(34)}\X_3 \cup P_{(34)}\X_j, P_{(34)}\X_3 \cup \X_j$ and $P_{(34)}\X_3 \cup P_{(34)}\X_j \cup \X_3 \cup \X_j,j=1,\ldots,4$ are groups with respect to matrix multiplication.


Now let $G$ represent any matrix group from the chain of groups corresponding to $\X_k.$ Now by Theorem \ref{thm1} we observe that if $D \in \X_k$ then there exist $B\in \Y_k$ and $C \in \Z_k$ such that $B=P_{(23)}DP_{(23)}$ and  $C = P_{(24)}DP_{(24)},$  $k=1,\ldots,4.$ 
Thus the chains corresponding to $\Y_3$ and $\Z_3$ follows from the observation that $f: G\rightarrow G',$ is a group isomorphism defined by  $f(M)=P_{(23)}MP_{(23)}$ and $f(M)=P_{(24)}MP_{(24)}$ when  $G'$ is the image of the map $f.$ Note that $P_{(23)}^T=P_{(23)},$ $P_{(24)}^T=P_{(24)}.$ $\hfill{\square}$

Then we have the following observations from Theorem \ref{prop:complex gp}.
\begin{remark}
\begin{enumerate}
 \item The groups of OPMs, $P_{(34)}\X_3,$ $P_{(24)}\Y_3$ and $P_{(23)}\Z_3$ 
 are commutative and a matrix $A$ in any of these groups can be written as $A=x I +y P+z P^2 +w P^3,$ where $P= P_{(1324)}$ if $A\in P_{(34)}\X_3,$ $P= P_{(1234)}$ if $A\in P_{(24)}\Y_3,$ and $P=P_{(1342)}$ if $A\in P_{(23)}\Z_3.$

\item $\O\P_4$ is not closed under matrix multiplication and hence $\O\P_4$ does not form a group: Consider $$A=\bmatrix{\frac{2}{5}& -\frac{2}{5}&\frac{4}{5}&\frac{1}{5}\\ -\frac{2}{5}& \frac{2}{5}& \frac{1}{5}& \frac{4}{5}\\ \frac{4}{5}& \frac{1}{5}& -\frac{2}{5}& \frac{2}{5}\\ \frac{1}{5} & \frac{4}{5} &\frac{2}{5}& -\frac{2}{5}} \in \X_1 \,\,\mbox{and}\,\, B=\bmatrix{\frac{\sqrt2}{3}& \frac{2}{3}& -\frac{\sqrt2}{3}&\frac{1}{3}\\\frac{2}{3}& -\frac{\sqrt2}{3}& \frac{1}{3}& \frac{\sqrt2}{3}\\ -\frac{\sqrt2}{3}& \frac{1}{3}&  \frac{\sqrt2}{3}& \frac{2}{3}\\  \frac{1}{3}& \frac{\sqrt2}{3}& \frac{2}{3}& -\frac{\sqrt2}{3}}\in \Y_1,$$
    then clearly $AB=\bmatrix{ -\frac{2\sqrt2}{15}-\frac{1}{5}& \frac{8}{15}+\frac{\sqrt 2}{5}& \frac{2\sqrt 2}{15}&  \frac{2}{3}-\frac{\sqrt2}{5}\\ -\frac{\sqrt2}{5}+\frac{8}{15}&  -\frac{1}{5}+\frac{2\sqrt2}{15}&\frac{\sqrt2}{5}+\frac{2}{3}&  -\frac{2\sqrt2}{15}\\ \frac{2\sqrt2}{5}+\frac{4}{15}& \frac{2}{5}+\frac{\sqrt2}{15}&-\frac{2\sqrt2}{5}+\frac{1}{3}& -\frac{\sqrt2}{15}\\ -\frac{\sqrt2}{15}+\frac{2}{5}& \frac{4}{15}-\frac{2\sqrt 2}{5}& \frac{\sqrt2}{15}& \frac{1}{3}+\frac{2 \sqrt2}{5}}\not\in \O\P_4.$ 


\end{enumerate}
\end{remark}

Then the following corollary describes all real OPMs.

\begin{corollary}\label{cor:real gp} (Characterization of real OPMs)
Under the assumptions and notations of Theorem \ref{thm1}, a matrix $A\in \X\cup \Y \cup \Z$ where \beano \X &=& \left\{ {\overline{P}} M^{s}_{x,z}, {\overline{P}}N^{s}_{z,x}: s=\pm\right\}\\ 
\Y &=& \left\{ {\overline{P}}P_{(23)}M^{s}_{x,z}P_{(23)}, {\overline{P}} P_{(23)}N^{s}_{z,x}P_{(23)} : s=\pm \right\}\\
\Z &=&  \left\{ {\overline{P}}P_{(24)}M^{s}_{x,z}P_{(24)}, {\overline{P}}P_{(24)}N^{s}_{z,x}P_{(24)} : s=\pm\right\} \eeano is a real OPM if and only if $x=\pm\sqrt{z(1-z)}, 0\leq z\leq 1$ for $s=+,$ and $x=\pm\sqrt{-z(1+z)}, -1\leq z\leq 0$  for $s=-.$

\end{corollary}

Observe that the parametric curves which define the real OPMs as given by Corollary \ref{cor:real gp}  are $x^2+z^2+ rz=0,$ $r\in\{1,-1\}.$ Then, $$(x, z)= \left(\frac{1}{2}\sin\theta, -\frac{r}{2}(1-r\cos\theta) \right), \,\, -\pi \leq \theta\leq \pi$$ provides one-parameter trigonometric parametrizations for the parametric curves. In particular, from equation (\ref{eqn:groverxyz}), trigonometric parametrizations of the continuous deformations of the Grover matrix of order $4$ can be obtained from the trigonometric parametrizations of sets of OPMs $\X_1, \Y_1, \Z_1$ given by
\be \label{coin:x1} {(\X_1)_{\theta}=\left\{\bmatrix{\frac{1}{2}\sin\theta&-\frac{1}{2}\sin\theta&\frac{1}{2}(1+\cos\theta)&\frac{1}{2}(1-\cos\theta)\\-\frac{1}{2}\sin\theta&\frac{1}{2}\sin\theta&\frac{1}{2}(1-\cos\theta)&\frac{1}{2}(1+\cos\theta)\\\frac{1}{2}(1+\cos\theta)&\frac{1}{2}(1-\cos\theta)&-\frac{1}{2}\sin\theta&\frac{1}{2}\sin\theta\\\frac{1}{2}(1-\cos\theta)&\frac{1}{2}(1+\cos\theta)&\frac{1}{2}\sin\theta&-\frac{1}{2}\sin\theta}:\theta \in [-\pi,\pi]\right\},}\ee 
\be \label{coin:y1}(\mathcal{Y}_1)_{\theta}=\left\{ \bmatrix{\frac{1}{2}\sin\theta &\frac{1}{2}(1+\cos\theta)&-\frac{1}{2}\sin\theta&\frac{1}{2}(1-\cos\theta)\\\frac{1}{2}(1+\cos\theta)&-\frac{1}{2}\sin\theta&\frac{1}{2}(1-\cos\theta)&\frac{1}{2}\sin\theta\\-\frac{1}{2}\sin\theta&\frac{1}{2}(1-\cos\theta)&\frac{1}{2}\sin\theta&\frac{1}{2}(1+\cos\theta)\\\frac{1}{2}(1-\cos\theta)&\frac{1}{2}\sin\theta&\frac{1}{2}(1+\cos\theta)&-\frac{1}{2}\sin\theta}:\theta \in [-\pi,\pi] \right\},\ee 
\be  \label{coin:z1}(\Z_1)_{\theta}=\left\{\bmatrix{\frac{1}{2}\sin\theta&\frac{1}{2}(1+\cos\theta)&\frac{1}{2}(1-\cos\theta)&-\frac{1}{2}\sin\theta\\\frac{1}{2}(1+\cos\theta)&-\frac{1}{2}\sin\theta&\frac{1}{2}\sin\theta&\frac{1}{2}(1-\cos\theta)\\\frac{1}{2}(1-\cos\theta)&\frac{1}{2}\sin\theta&-\frac{1}{2}\sin\theta&\frac{1}{2}(1+\cos\theta)\\-\frac{1}{2}\sin\theta&\frac{1}{2}(1-\cos\theta)&\frac{1}{2}(1+\cos\theta)&\frac{1}{2}\sin\theta}:\theta \in [-\pi,\pi]\right\}\ee respectively.

Next, in what follows, we characterize all rational OPMs. Treating $x^2+z^2-z=0$ as a polynomial in indeterminate $z$, we obtain
$$z=\frac{1\pm \sqrt{1-4x^2}}{2}.$$
Then $z \in \mathbb{Q}$ if and only if $1-4x^2$
 is zero or perfect square of a nonzero rational number of the form $p/q$ for $p,q \in \mathbb{Z}, q\neq 0.$ It is zero if $x \in \{-\frac{1}{2},\frac{1}{2}\}.$ If $1-4x^2=\frac{p^2}{q^2}$ for $p/q \neq 0$ then after rewriting it takes  the form $X^2-4Y^2=1,$ where $X=q/p$ and $Y=x q/p.$ Now $(X+2Y)$ and $(X-2Y)$ are units in $\mathbb{Q}$ for $x \in \mathbb{Q}$ such that $(X+2Y)(X-2Y)=1.$ Thus letting $X+2Y=r$ and $X-2Y=\frac{1}{r}$ for some nonzero $r \in \mathbb{Q},$ we obtain  
 $$X=\frac{1}{2}\left(r+\frac{1}{r}\right)~\mathrm{and}~Y=\frac{1}{4}\left(r-\frac{1}{r}\right),$$ which ultimately gives us the values of $x$ and $z$ in terms of the parameter $r.$ Similar procedure can be followed for  $x^2+z^2+z=0.$ Thus we have the following corollary. 

 \begin{corollary}\label{cor:rational gp} (Characterization of rational OPMs)
Under the assumptions and notations of Corollary \ref{cor:real gp}, a matrix $A\in \mathcal{X}\cup\mathcal{Y}\cup\mathcal{Z}$ is a rational  OPM if and only if $x=\dfrac{r^2-1}{2(r^2+1)},z=\dfrac{1}{2}\pm \dfrac{r}{r^2+1}$ for $s=+$ and $x=\dfrac{r^2-1}{2(r^2+1)},z=-\dfrac{1}{2}\pm \dfrac{r}{r^2+1},$ for $s=-$, where $r\in \Q.$

 \end{corollary}
 
 Finally, we mention that the chains of matrix groups described in Theorem \ref{prop:complex gp} remains valid when the subgroups are restricted only for the real or rational matrices as given in Corollary \ref{cor:real gp} and  Corollary \ref{cor:rational gp}.

\section{Search for orthogonal matrices that are linear combinations of permutation matrices but not permutative}

In \cite{sarkar2020}, it is shown that an orthogonal matrix of order $3$ is a linear combination of permutation matrices if and only if it is a permutative matrix. However, this is no longer true for matrices of order $4$ as follows from the following example. Consider the block diagonal matrix $$A=
 \left[
\begin{array}{cc}
     1 & \mathbf{0} \\
    \mathbf{0}^T & \frac{2}{3} J_3-I_3 \\
   \end{array}
   \right]
   =-\frac{1}{3} I+\frac{2}{3}P_{(234)}+\frac{2}{3}P_{(243)},$$ where $J_3=\boldsymbol{1}_3\boldsymbol{1}_3^T.$ Then $A$ is an orthogonal matrix which is a linear combination of permutation matrices but not a permutative matrix. Indeed, it may be noted that this matrix $A$ is a direct sum of two permutative matrices. Then the following question arises: Does there exist an orthogonal matrix of order $4$ which is a linear combination of permutation matrices but neither a permutative matrix nor a direct sum of permutative matrices with permutations of its rows and columns? In this section, we investigate this problem.
   

First we derive certain sufficient conditions for which an orthogonal matrix which is a linear combination of permutation matrices is always permutative, that is, such a matrix can always be written in the form given by equation (\ref{PO}). Also, recall that a necessary condition for a linear combination of permutation matrices to be real orthogonal is that sum of the entries along a row or column should be $\pm1$ \cite{Kapoor1975}. We provide an alternative easy proof of this result for orthogonal matrices of order $4$ in the following proposition. 

\begin{proposition} \label{proposition:necessary condition}
 A necessary condition for a linear combination of permutation matrices $A$ of order $4$ to be orthogonal is that the sum of the entries of $A$ along each row and column is $\pm 1.$
 \end{proposition}
 
 \pf Suppose $A=\sum_{\sigma\in S_4} x_{\sigma}P_{\sigma},$ where $S_4$ denotes the symmetric group of order $4.$ Then the $i$th row sum of $A$ is
 $$A_{(i,:)} = \sum_{j=1}^{4}\sum_{\sigma}{x_{\sigma}P_{\sigma}(i,j)}
   = \sum_{\sigma}{x_{\sigma}}\sum_{j=1}^{n}P_{\sigma}(i,j)
   =\sum_{\sigma} {x_{\sigma}}.$$
Similarly, the $i$th column sum $A_{(:,i)}$ of $A$ is $\sum_{\sigma} {x_{\sigma}}.$
Then consider the Hardamard matrix of order $4$ as follows: \begin{equation}\label{hadamard}
    H=\frac{1}{2}\bmatrix{1&1&1&1\\1&-1&1&-1\\1&1&-1&-1\\1&-1&-1&1}.
\end{equation}

Setting $ B=H A H,$ $Be_1=(H A H)e_1=(H A)(He_1)=\frac{1}{2}H(A\boldsymbol{1})=\frac{1}{2}\sum_{\sigma} {x_{\sigma}}(H\boldsymbol{1})=(\sum_{\sigma} {x_{\sigma}})e_1,$ where $e_1=\bmatrix{1&0&0&0}^{T}$ and $\boldsymbol{1}=\bmatrix{1&1&1&1}^{T}.$ Similarly, 
$B^{T}e_1=(\sum_{\sigma}{x_{\sigma}})e_1.$ 
Then $$
 B=HAH =\bmatrix{\sum_{\sigma} {x_{\sigma}}&\mathbf{0}\\
 \mathbf{0}^T&\bar{A}},
 $$
 where $\mathbf{0}=\bmatrix{0&0&0}$ and $\bar{A}$ is a $3 \times 3$ orthogonal matrix. Consequently,
 $\sum_{\sigma} {x_{\sigma}}\in\{\pm 1\}$ since $B$ is orthogonal. This completes the proof.
 $\hfill{\square}$


We call a set of real (nonzero) matrices of order $k$, $S=\{A_1, A_2, \hdots, A_n\}$ is pairwise $H$-orthogonal if Hadamard product of any pair of matrices $A_i, A_j, i\neq j$, denoted by $A_i\circ A_j$ is the zero matrix, $1\leq i,j\leq n.$ We denote $\langle S\rangle = \left\{\sum_{j=1}^n \alpha_jA_j : \alpha_j\in \R\right\}$ as the vector space generated by elements of $S.$ Observe that if $S$ is a pairwise $H$-orthogonal set of permutation matrices then $A\in \langle S \rangle$ is always permutative.



Then we have the following proposition. 

\begin{proposition} \label{lemma:6 Lsum}
Any linear combination of permutation matrices of order $4$ can be written as a  sum of at most $6$ permutative matrices each of which is a linear combination of $4$ permutation matrices that are pairwise $H$-orthogonal.  
\end{proposition}
\pf The proof follows from the partition of the symmetric group $S_4=\cup_{k=1}^6 \tilde{S}_k$ where  \beano  &&  \tilde{S}_1=\{id,(12)(34),(13)(24),(14)(23)\}, \,\,  \tilde{S}_2=\{(23),(124),(1342),(143)\}, \\
&& \tilde{S}_3 =\{(24),(123),(134),(1432)\}, \,\, \tilde{S}_4=\{(34),(12),(1324),(1423)\}, \\
&& \tilde{S}_5 = \{(14),(1243),(132),(234)\}, \,\,
\tilde{S}_6= \{(13),(1234),(142),(243)\}.
\eeano such that $M_{\tilde{S}_k}=\{P_\sig : \sig\in \tilde{S}_k\}$ is a pairwise $H$-orthogonal set and any $A\in \langle M_{\tilde{S}_k} \rangle$ is a permutative matrix. \hfill{$\square$}




Then we have the following theorem.

\begin{theorem} \label{proposition:LS+P}
Let $A\in\langle S\rangle$ where $S$ is a pairwise $H$-orthogonal set of permutation matrices of order $4.$ Then $A+c P$  is an OPM for any $c\in \R$ and $P\in \P_4.$  
 \end{theorem}
 \pf Obviously $B=A\in \langle S\rangle$ is permutative if $c=0$. Let $c\neq 0.$ Let $A$ be a symbolic permutative matrix with first row ${\bf x}= (x, y, z, w).$ Consider the entries $A_{ij}$ for which $P_{ij}=1.$ Then $B_{ij}=A_{ij}+c$ if and only if $P_{ij}=1;$ and $B_{ij}=A_{ij}$ otherwise. For any pair of indices $(i,j)$ and $(k,l)$ with $P_{ij}=1=P_{kl},$ the unit norm condition of rows of $B$ implies $A_{ij}=A_{kl}$ since $c\neq 0.$ Thus the permutative structure of $A$ implies that $B$ is permutative.  $\hfill{\square}$

Below we show that no orthogonal matrix of order $4$ can be a linear combination of two  distinct permutation matrices. By non-trivial linear combination we mean all the coefficients of the linear combination are non-zero.

\begin{theorem} \label{Proposition:lsum 2permutations}
There is no orthogonal matrix which is a non-trivial linear combination of two distinct permutation matrices.
 \end{theorem}
 \pf Let $A=\alpha P + \beta Q$ be an orthogonal matrix and $P\neq Q.$ If $P \circ Q =0$ then $A$ is an OPM. From the classification of all OPMs of order $4$ described in Remark \ref{remark1}, any OPM is a linear combination of four $H$-orthogonal permutation matrices. 
 Indeed, it follows from equations (\ref{eqn:variety x}), (\ref{eqn:variety y}), (\ref{eqn:variety z}) that if one or more coefficients in the linear combination of $H$-orthogonal permutation matrices are zero, then the corresponding OPM becomes $\pm 1$ times a permutation matrix. Hence the desired result follows.
  
 Next, assume that $P\circ Q\neq 0$ and $A=\alpha P + \beta Q$ is an orthogonal matrix, $P,Q\in\P_4.$ This means there can exist at most two pairs of indices $(i,j)$ such that $P_{ij}=Q_{ij}=1$ since $P\neq Q.$ Then $A_{ij}=\alpha+\beta$ for those $(i,j),$ and two permutation matrices $X,Y$ can be found for which $XAY=(\alpha+\beta)I_1 \oplus A_1$ or $(\alpha+\beta)I_2 \oplus A_2$ where $A_1$ and $A_2$ are OPMs of orders $3$ and $2$ respectively. Indeed, each row of $A_1$ is a permutation of $\{0, \alpha, \beta\};$ whereas each row of $A_2$ is a permutation of $\{\alpha, \beta\}.$ Then from the classification of OPMs of order $3$ (see Theorem 3.1, \cite{sarkar2020}) it can be seen that either $(\alpha, \beta)=(\pm 1,0)$ or $(\alpha,\beta)=(0,\pm 1).$ The same holds for $A_2,$ and hence the desired result follows. $\hfill{\square}$  
 
  The following theorem provides characterization of orthogonal matrices that are linear combinations of three permutation matrices.

\begin{theorem} \label{Proposition:lsum 3permutations}
If an orthogonal matrix $A$ is a (real) linear combination of three distinct  permutation matrices then either $\pm A$ is a permutation matrix or $XAY$ is a direct sum of OPMs of order $3$ and $1$, for some permutation matrices $X,Y.$
\end{theorem}
 \pf Let $A=\alpha P + \beta Q + \gamma R$ be an orthognal matrix, where $P, Q, R$ are distinct permutation matrices of order $4$. Then two cases arise. Either the set $S=\{ P, Q,R\}$ is pairwise $H$-orthogonal or there exist at least one pair of matrices in $S$ whose Hadamard product is a non-zero matrix.  Note that sum of entries of each row and column of $A$ is $\alpha+\beta +\gamma.$
 
 First suppose that $S$ is a pairwise $H$-orthogonal set. Then following a similar argument as in the proof of Theorem \ref{Proposition:lsum 2permutations} it can be  concluded that $A=\pm M$ for some  $M \in \P_4.$ Next assume that $S$ is not pairwise $H$-orthogonal. If there is only one pair of elements of $S$ that are not $H$-orthogonal. Without loss of generality, let $P\circ Q\neq0.$ Then $A_{ij}=\alpha+\beta$ if and only if $P_{ij}=Q_{ij}=1$ for at most two indices $(i,j).$ Hence the unit norm condition of rows and columns of $A$ yields the polynomial system:
 $$ {\alpha}^2+{\beta}^2+{\gamma}^2=1, \,\,\, ({\alpha}+{\beta})^2+{\gamma}^2=1.$$
  Then it is computational to check that either ${\alpha}=0$ or ${\beta}=0.$ Then from Theorem \ref{Proposition:lsum 2permutations} it follows that $A=\pm M$ for some $M\in \P_4.$ Now assume that there are two distinct pairs of matrices in  $S$ each of which are not $H$-orthogonal. Without loss of generality, let  $P\circ Q\neq 0,P\circ R\neq 0.$ Then $A_{ij}=\alpha+\beta$ and $A_{kl}=\alpha+\gamma$ if and only if $P_{ij}=Q_{ij}=1,$ $P_{kl}=R_{kl}=1$ for some indices $(i,j)$ and $(k,l).$ Thus  
    $(\alpha,\beta,\gamma)$ satisfy the following polynomial system due to the unit norm condition of rows of $A$: $$(\alpha+\beta)^2+\gamma^2=1, \,\,\,  (\alpha+\gamma)^2+\beta^2=1. $$ Solving these equations we have either
 $\alpha=0$ or $\beta=\gamma.$ If $\alpha=0,$ $A$ is linear combination of $Q$ and $R$ and hence the desired result follows from Theorem \ref{Proposition:lsum 2permutations}, while for $\beta=\gamma$ rows of $A$ are permutations of $\mathbf{x}_1=(\alpha+\beta,\beta,0,0)$ or $\mathbf{x}_2=(\alpha,\beta,\beta,0).$ If rows of $A$ are permutations of $\mathbf{x}_1$ only then $\beta({\alpha}+{\beta})=0.$ Otherwise if permutations of both $\mathbf{x}_1$ and $\mathbf{x}_2$ present as rows of $A$ we obtain ${\alpha}{\beta}=0.$  Thus the result follows from Theorem \ref{Proposition:lsum 2permutations}. Finally, let all pairs of matrices from $S$ are not $H$-orthogonal. Then if $P\circ Q \circ R \neq 0$ then there is exactly one index $(i,j)$ such that $A_{ij}=\alpha+\beta+\gamma$ and $P_{ij}=Q_{ij}=R_{ij}=1$ since $P\neq Q \neq R,$ which further implies that $\alpha+\beta+\gamma \in \{\pm1\}.$ Obviously, two permutation matrices $X,Y$ can be found such that $XAY=(\alpha+\beta+\gamma) \oplus A_1$ where $A_1$ is an orthogonal matrix of order $3$ which is linear combination of $3$ permutation matrices. Then using the characterization of orthogonal matrices that are linear combinations of permutations (see Theorem $3.2$, \cite{sarkar2020}) we conclude that $A_1$ is a OPM and the desired result follows. Otherwise, if  $P \circ Q \neq {0},Q \circ R \neq {0},R \circ P \neq {0}$ with $P\circ Q \circ R=0$ then each rows of $A$ can be permutations of $(\alpha+\beta,\gamma,0,0)$, $(\alpha+\gamma,\beta,0,0)$ and $(\alpha,\beta+\gamma,0,0).$ However, all of these row vectors can not appear as rows of $A$ simultaneously since column sum of $A$ is $\alpha+\beta+\gamma$ for each column.  These complete the proof.
 $\hfill{\square}$

Now we focus on finding orthogonal matrices  that are real linear combinations of permutations but neither permutative nor a direct sum of permutative matrices up to permutations of rows or columns. Thus we investigate existence of such matrices belonging to subspaces and direct sum of subspaces of matrices of order $4$ that are generated by sets of linearly independent permutation matrices. Recall that the real linear combinations of permutation matrices of order $n$ form a vector space of dimension $(n-1)^2+1$ \cite{farahat1970sets}. We denote this space as $\L$ for $n=4.$ We choose a basis $\mathcal{B}$ of $\L$ that contains $10$ permutation matrices given by 
$$\mathcal{B}=\{P_{(12)}, P_{(23)}, P_{(24)}, P_{(34)}, P_{(123)}, P_{(124)}, P_{(234)}, P_{(12)(34)}, P_{(13)(24)}, P_{(14)(23)}\}.$$

Let $A \in \L.$ Then
$A \in \oplus_{k=1}^{5} \L_k,$ $\L_k=\langle \B_k\rangle$ is the subspace generated by $\B_k, k=1,\ldots,5$ as follows: \beano
&& \B_1=\{P_{(12)},P_{(34)},P_{(13)(24)},P_{(14)(23)}\},
 \B_2=\{P_{(24)},P_{(12)(34)}\},\\
&& \B_3=\{P_{(124)},P_{(234)}\},
 \B_4=\{P_{(123)}\},
 \B_5=\{P_{(23)}\}.
 \eeano Note that $\B_k$ is a pairwise $H$-orthogonal set for each $k.$ In particular, if $A\in\L_k, k=1,\ldots,5,$ is orthogonal then $A$ can be characterized by Theorem \ref{thm1} and hence $A \in \O\P_4.$
 

 
 
 Now we briefly review the concept of \textit{combinatorially orthogonal} matrices introduced by Brualdi et al. \cite{Brualdi1991} which will be used in sequel. A matrix having entries from $\{0,1\}$ is called a $(0,1)~matrix$. The nonzero \textit{pattern} of a matrix $A$ is defined as a $(0,1)~ matrix$ $M_A$ such that $ij$th entry of $M_A=1$ if and only if $a_{ij}\neq 0$. A nonzero \textit{pattern} $M$ is orthogonal if there exists a (real) orthogonal matrix with the same pattern.  Let $A$ be a $(0,1) ~matrix$ of order $n.$ Then $A$ is \textit{combinatorially orthogonal} or \textit{quadrangular} if inner product of distinct rows or columns is not equal to $1.$ Let $S$ be a subset of rows of $A$ such that for each element of $S$ there is  another element of $S$ with nonzero inner product. Then $A$ is said to be row strongly quadrangular if the matrix, whose rows are all the elements of $S$, has at least $|S|$ number of columns with at least two $1$s. Similarly, the matrix $A$ is said to be column strongly quadrangular if the set $S$ contains columns of $A$ and if the matrix whose columns are all the elements of $S$, has at least $|S|$ number of rows with at least two $1$'s. If a $(0,1)$ matrix is both row and  column strongly quadrangular then it is called \textit{strongly quadrangular}. 
 
 Note that if a $(0,1)$ matrix supports unitary then it is strongly quadrangular but the converse need not be true. 
 Now we recall the following proposition from \cite{severini2008}.
 
 \begin{proposition} \label{prop:st qua} A $(0,1)$ matrix of degree $n \leq 4 $ supports a unitary if and only if it is strongly quadrangular.\end{proposition}
 
 Then we have the following theorem.
 
  \begin{theorem} \label{theorem:L_i+L_j}
  Let $A \in \L_i \oplus \L_j,$ $i,j\in\{1,\ldots,5\}$ be an orthogonal matrix.   Then $A \in \O\P_4.$ 
  \end{theorem}
 \pf It is clear from Theorem \ref{proposition:LS+P} and Theorem \ref{thm1} that $A \in \O\P_4$ whenever $A \in \L_i \oplus \L_j$ for $i\neq j,i\in\{1,\ldots,5\},$ $j\in\{4,5\}.$
 
 First let $A=a_1P_{(12)}+a_2P_{(34)}+a_3P_{(13)(24)}+a_4P_{(14)(23)}+b_1 P_{(24)}+b_2 P_{(12)(34)} \in \L_1 \oplus \L_2.$
 Then the unit norm conditions of $2$nd and $4$th rows, $1$st and $2$nd rows, and $3$rd and $4$th rows of $A$ yield $b_2=0$ or $a_1=a_2,$ $b_1=0$ or $a_2=a_3$, and $b_1=0$ or $a_1=a_3,$ respectively.
 If $b_1=b_2=0,$ then $A\in \L_1.$ If $b_1\neq 0,b_2=0$ and $a_1=a_2=a_3$ then $A=A\left(\mathbf{x};P_{(1432)},P_{(13)(24)},P_{(1234)}\right),$ where $\mathbf{x}=(a_1+b_1,a_1,a_1,a_4).$ If $b_1=0,b_2\neq 0,a_2=a_1,$ then $A=A\left(\mathbf{x};P_{(12)(34)},P_{(13)(24)},P_{(14)(23)}\right),$ where $\mathbf{x}=(a_1,a_1+b_2,a_3,a_4).$ 
At last, for $b_1,b_2\neq 0, a_1=a_2=a_3,$ $A=A\left(\mathbf{x};P_{(1432)},P_{(13)(24)},P_{(1234)}\right),$ where $\mathbf{x}=(a_1+b_1,a_1+b_2,a_1,a_4).$
Thus by Theorem \ref{thm1} in all the above cases $A\in \O\P_4.$ 

If $A\in \L_1 \oplus \L_3$ then $A$ is linear combination of at most $6$ permutations and using similar arguments as above it is easy to verify that $A\in \O\P_4.$
     
     Next, suppose $A=b_1P_{(24)}+b_2P_{(12)(34)}+c_1P_{(124)}+c_2P_{(234)}\in \L_2 \oplus \L_3.$ Then the $(0,1)$ pattern $M_A$ of $A$ is given by
     $$M_{A}=\bmatrix{1&1&0&0\\1&0&1&1\\0&0&1&1\\1&1&1&0},$$ which is not quadrangular. Hence $A$ is not an orthogonal matrix with the $(0,1)$ pattern $M_A$. Let $M_i$ denote the $i$th column of  $M_A.$ However, for  $A$ to be orthogonal  which means $M_A$ to be quadrangular, the coefficients must satisfy the following conditions: $b_2=0$ or $b_1+c_1=0$, and $b_2=0$ or $b_1+c_2=0$ by setting $M_1^TM_4=0$ and $M_2^TM_3=0$ respectively.

     If $b_2=0,$ then by Proposition \ref{prop:st qua} the nonzero pattern of $A$ can support an orthogonal matrix. Further unit norm condition of rows of $A$ implies either $b_1=0$ or $c_1=c_2.$ If $b_2=0$ and $b_1=0,$  $A \in \L_3.$ Otherwise for $c_1=c_2$ together with $b_2=0,$ $A$ is a permutative matrix with two nonzero entries in each row. Hence from Theorem \ref{Proposition:lsum 2permutations}, $A=\pm R,R\in \P_4.$
         If $b_2\neq 0$ i.e. $b_1+c_1=0$ and $b_1+c_2=0,$ a further analysis yields $c_1=c_2=0.$ Hence $A=b_2 P_{(12)(34)},$ where $b_2=\pm 1.$
      These complete the proof. $\hfill{\square}$


\begin{theorem} \label{thm1:3L}
   Let $A \in \L_1\oplus \L_i \oplus \L_j$ be orthogonal where $i,j\in\{2,3,4,5\}$ and $(i,j)\notin \{(2,5),(3,4)\}.$ Then $A\in\O\P_4.$
 \end{theorem}
 
 \pf Since $ \L_1 $ is generated by four $H$-orthogonal permutation matrices and each of $\L_i$ and $\L_j$ is generated by one or two permutation $H$-orthogonal matrices, the nonzero pattern of $A$ is the all-one matrix in general. 
 Let $A=a_1P_{(12)}+a_2P_{(34)}+a_3P_{(13)(24)}+a_4P_{(14)(23)}+b_1P_{(24)}+b_2P_{(12)(34)}+c_1P_{(124)}+c_2P_{(234)}\in \L_1\oplus \L_2 \oplus \L_3.$ Then the unit norm conditions of $1$st row and $1$st column, $3$rd row and $3$rd column, $3$rd row and $4$th column, and $1$st row and $2$nd column of $A$ yield $c_1=0$ or $a_4=a_1+b_2,$ $c_2=0$ or $a_4=a_2+b_2,$ $b_1+c_1=0$ or $a_3=a_1,$ and $b_1+c_2=0$ or $a_3=a_2,$ respectively.  Then it can be verified  that: 
 \begin{enumerate}
     \item If $c_1=c_2=0,$ then $A\in \L_1 \oplus \L_2.$
     \item Consider $c_1=0,c_2 \neq 0,a_4=a_2+b_2,b_1+c_1=0,a_3=a_2.$
     Further, $b_2=0$ or $a_1=a_2$ and $c_2=0$ or $a_2=a_4,$ and hence when $b_2=0,$ $A \in \L_1 \oplus \L_3.$ Otherwise we have $a_1=a_2.$ Hence $a_1=a_2=a_4,$ which again implies $b_2=0.$ Thus $A\in \L_1 \oplus \L_3.$
     \item If $c_1=0,c_2 \neq 0,a_4=a_2+b_2,b_1+c_2=0,a_3=a_1,$ then we obtain  $b_2=0$ or $a_1=a_2.$  Further it can be verified that in both the cases $A$ can not be an orthogonal matrix under all the given conditions.
     
     \item When $c_1=0,c_2 \neq 0,a_4=a_2+b_2,a_3=a_2=a_1,$ then the orthonormality of $A$ further implies $a_4=a_1+b_1$ and hence $b_1=b_2.$ So that $A$ becomes $A\left(\mathbf{x};P_{(1324)},P_{(1423)},P_{(12)(34)}\right)$ with $\mathbf{x}=(a_4+c_2,a_4,a_1,a_4).$
     A similar analysis can be done for all the cases where $c_2=0$ and $c_1 \neq 0.$
     \item If $c_1\neq 0,c_2 \neq 0,a_4=a_1+b_2,a_1=a_2,b_1+c_1=0,b_1+c_2=0,$ then $A$ becomes $A\left(\mathbf{x};P_{(1234)},P_{(13)(24)},P_{(1423)}\right)$ with $\mathbf{x}=(a_1,a_4+c_1,a_3,a_4).$ 
     \item Consider $c_1\neq 0,c_2 \neq 0, a_4=a_1+b_2, a_1=a_2=a_3, b_1+c_1=0.$ Then the orthogonality condition of $A$ yields $a_1=0$ or $2a_4+c_2=0.$ For $a_1=0$ we obtain $a_4=0.$ Hence $A \in \L_2 \oplus \L_3$ or $b_1+c_2=0,$ which satisfies the previous case. Otherwise for $2a_4+c_2=0$ we obtain  $b_1=2a_4$, which is equivalent to say $b_1+c_2=0,$ or $b_1=0,$ hence it follows from Proposition \ref{proposition:necessary condition} that $A=\pm(\frac{1}{2} J_4-P_{(234)}).$ Similarly it can be done for $c_1\neq 0,c_2 \neq 0,a_4=a_1+b_2,a_1=a_2=a_3,b_1+c_2=0.$
     
     \item If $c_1\neq 0,c_2 \neq 0,a_4=a_1+b_2,a_1=a_2=a_3,$ then the orthonormality of $A$ yields $c_1=c_2$ or $a_1+b_1=a_4.$ Now, $c_1=c_2$ implies $A=A\left(\mathbf{x};P_{(14)(23)},P_{(13)(24)},P_{(12)(34)}\right)$ with $\mathbf{x}=(a_1+b_1+c_1,a_4+c_1,a_1,a_4).$ In particular $A=\pm(\frac{1}{2}J_4-P),P\in \P_4.$ On the other hand for $a_1+b_1=a_4$ we have $b_1=b_2$ and $A$ takes the form $A\left(\mathbf{x};P_{(1324)},P_{(1423)},P_{(12)(34)}\right)$ with $\mathbf{x}=(a_4+c_2,a_4+c_1,a_1,a_4).$ 
 \end{enumerate}
 
Using similar arguments of step by step elimination of entries of the polynomial system defined by orthonormality condition of the columns and rows of concerned matrices it can be verified that any orthogonal matrix belongs to
$\L_1 \oplus \L_2 \oplus \L_4, \L_1 \oplus \L_3 \oplus \L_5,$ or $\L_1 \oplus \L_4 \oplus \L_5$ is in $\O\P_4.$  $\hfill{\square}$
 
  \begin{theorem} \label{thm2:3L} 
   Let $A \in \L_i\oplus \L_j\oplus\L_k$ be orthogonal where $i,j,k\in\{2,3,4,5\},$ then $A\in\O\P_4.$
 \end{theorem}
 \pf Clearly for all the given choices of $i,j,k,$ $A$ can be a linear combination of at most $5$ permutation matrices and nonzero patterns of $A$ are  other than the all-one  matrix. 
 Suppose, $A=b_1P_{(24)}+b_2P_{(12)(34)}+c_1P_{(124)}+c_2P_{(234)}+eP_{(23)}\in \L_2 \oplus \L_3\oplus\L_5.$ Clearly the $(0,1)$ pattern of $A$ is not quadrangular and $b_2=0$ or $b_1+c_2+e=0$ and $e=0$ or $b_2+c_1=0$ are satisfied. Hence if $e=0$ then $A \in \L_2 \oplus \L_3.$
 Otherwise if $e\neq 0$ and $b_2=0$ then $c_1=0,$ so that nonzero pattern of $A,$
$$M_{A}=\bmatrix{1&0&0&0\\0&0&1&1\\0&1&1&1\\0&1&0&1}$$ is non-quadrangular. Then considering the 3rd and 4th rows of $A$ to be orthogonal to each other, it follows that $b_1$ should be $0.$  Thus $A\in \L_3 \oplus \L_5.$
 Finally, $b_1+c_2+e=0$ and $b_2+c_1=0$ can not hold simultaneously  since $b_1+b_2+c_1+c_2+e\in \{\pm1\}$ is a necessary condition for $A$ to be orthogonal. Hence the proof follows from Theorem \ref{theorem:L_i+L_j}.
 
Similarly by looking into the nonzero patterns and eliminating some entries for the requirement of orthogonality of $A$ the desired result can be proved when $A$ belongs to each of the spaces $\L_2 \oplus \L_3\oplus\L_4,\L_2 \oplus \L_4\oplus\L_5,\L_3 \oplus \L_4\oplus\L_5.$ $\hfill{\square}$

 \begin{theorem}  Let $A \in \L_2\oplus\L_3\oplus\L_4 \oplus \L_5$ be orthogonal. Then $A \in \O\P_4.$ 
 \end {theorem}
 \pf Let $A=b_1P_{(24)}+b_2P_{(12)(34)}+c_1P_{(124)}+c_2P_{(234)}+d P_{(123)}+e_{(23)} P_{(23)} \in \L_2 \oplus \L_3\oplus\L_4 \oplus \L_5.$ Let $M_A^k$ for $k=1,\ldots,4$ denote the $(0,1)$ pattern of $A$ arise at different stages. Now,
$$M_{A}^1=\bmatrix{1&1&0&0\\1&0&1&1\\1&1&1&1\\1&1&1&1},$$ which is not quadrangular. Since  $ R_1R_2^T=1$ where $R_i$ denotes the $i$th row of $M_A^1,$ $b_2=0$ or $b_1+c_2+e=0$ holds. 

Hence if $b_2=0,$ then $$M_{A}^2=\bmatrix{1&1&0&0\\0&0&1&1\\1&1&1&1\\1&1&0&1},$$ which is again not quadrangular since $C_1^TC_3=1=C_2^TC_3,$ where $C_i$ denotes the $i$th column of $M_A^2.$ Thus $d=0$ or $b_1+c_1=0$ and $e=0$ or $b_1+c_1=0$ are satisfied.  If $d=e=0,$ then $A \in \L_2\oplus\L_3.$     If $d=0,e \neq 0$ then $A \in \L_2\oplus\L_3\oplus \L_5.$  Otherwise if $d\neq 0,e = 0$ then $A \in \L_2\oplus\L_3\oplus \L_4.$  Finally if $d\neq 0,e \neq 0,b_1+c_1=0,$ then          $$M_{A}^3=\bmatrix{1&1&0&0\\0&0&1&0\\1&1&0&1\\1&1&0&1}$$ and the orthogonality condition of corresponding $A$  implies $d=0$ or $e=0,$ which is a contradiction. Hence at least one of $d=0$ or $e=0$ holds whenever $b_2=0.$       

If $b_2 \neq 0$ then $b_1+c_2+e=0$ holds. Thus we obtain      $$M_{A}^4=\bmatrix{0&1&0&0\\1&0&1&1\\1&1&1&1\\1&1&1&1},$$ which is not quadrangular. Hence $e=0$ or $b_2+c_1+d=0$ and $b_1+c_2=0,$ or $b_2+c_1+d=0$ hold. So that $A\in \L_2\oplus\L_3\oplus \L_4$ if $e=0.$ While for $e\neq0$ we get $b_2+c_1+d=0$ and thus $b_1+b_2+c_1+c_2+d+e=0.$ Which leads to a contradiction by Proposition \ref{proposition:necessary condition}. Thus the proof follows from Theorem \ref{theorem:L_i+L_j} and Theorem \ref{thm2:3L}. 
$\hfill{\square}$

In all the cases above we determine spaces generated by specific permutation matrices such that any orthogonal matrix belonging to these spaces is either permutative or direct sum of permutative matrices after pre and post multiplication by permutation matrices to the original matrix. In the following we determine classes of orthogonal matrices that are linear combinations of permutation matrices but are not permutative matrices.

 \begin{theorem} \label{theorem:L_i+L_j+L_k}
   Let $A \in \L_1 \oplus \L_3 \oplus \L_4$ be orthogonal. Then either $A\in \O\P_4$ or there exist $P,Q\in \P_4$ such that $PAQ$ or $H(PAQ)H$ is of the form $\bmatrix{\pm 1&\mathbf{0}^T \\ \mathbf{0}&X}$ for some OPM $X$ of order $3,$  where $H$ is the Hadamard matrix of order $4$ as given in equation (\ref{hadamard}). 
 \end {theorem}
 
 \pf Suppose  $A=a_1P_{(12)}+a_2P_{(34)}+a_3P_{(13)(24)}+a_4P_{(14)(23)}+c_1P_{(124)}+c_2P_{(234)}+d P_{(123)} \in \L_1 \oplus \L_3 \oplus \L_4.$ Then the unit norm conditions of $1$st row and $2$nd column, $2$nd row and $3$rd column, $3$rd row and $1$st column of $A$ yield $c_2=0$ or $a_2=a_3$, $c_1=0$ or $a_1=a_3$ and $c_1=0$ or $a_1=a_4,$ respectively.
 Thus if $c_1 \neq 0$ and $c_2 \neq 0$ then we get $a_1=a_2=a_3=a_4$ and a further computation for $2$nd and $4$th rows gives $d=0,$ which implies $A \in \L_1\oplus \L_3.$ If $c_1=c_2=0,$ then $A\in \L_1 \oplus \L_4.$
If $c_1=0,a_2=a_3$ while $c_2 \neq 0$ then a further unit norm condition of $1$st and $2$nd columns lead to either $a_1=a_2$ or $d=0.$ If $d=0$ then we are done. Suppose $a_1=a_2,$ with $c_1=0,a_2=a_3,$ so that $a_1=a_2=a_3.$ Then orthogonality of $1$st and $2$nd rows of $A$ implies $a_1=0$ or $ a_1+ a_4+ c_2+d=0.$
    \begin{enumerate}
        \item Now while $a_1=0,$ we obtain
    $$P_{(12)}AP_{(13)}=\bmatrix{a_4+c_2+d&0&0&0\\0&d&c_2&a_4\\0&a_4&d&c_2\\0&c_2&a_4&d},$$ and hence $a_4+c_2+d=\pm 1.$ 
    \item For $a_1+ a_4+ c_2+d=0$ using Proposition \ref{proposition:necessary condition} we obtain $3 a_1+a_4+c_2+d=\pm 1$ so that 
    $2a_1=\pm 1$ and $a_4+c_2+d=-a_1.$ Hence  when row sum of $A$ is $1$ and $-1,$ we obtain $a_1=\frac{1}{2}$ and $a_1=-\frac{1}{2},$ respectively. Thus $P_{(12)}AP_{(13)}$ belongs to one of the following sets:
   \beano
   \mathcal{C}_1&=& \left\{\bmatrix{-\frac{1}{2}&\frac{1}{2}&\frac{1}{2}&\frac{1}{2}\\\frac{1}{2}&-a_4-c_2&\frac{1}{2}+c_2&a_4\\\frac{1}{2}&a_4&-a_4-c_2&\frac{1}{2}+c_2\\\frac{1}{2}&\frac{1}{2}+c_2&a_4&-a_4-c_2}: a_4=-\frac{1}{2}c_2\pm \frac{1}{2}\sqrt{(1-3 c_2)(1+c_2)}, \right. \\ && \hfill{\left. -1\leq c_2\leq \frac{1}{3}\right\}},\\
  \mathcal{C}_2&=& \left\{\bmatrix{\frac{1}{2}&-\frac{1}{2}&-\frac{1}{2}&-\frac{1}{2}\\-\frac{1}{2}&-a_4-c_2&-\frac{1}{2}+c_2&a_4\\-\frac{1}{2}&a_4&-a_4-c_2&-\frac{1}{2}+c_2\\-\frac{1}{2}&-\frac{1}{2}+c_2&a_4&-a_4-c_2}:a_4=-\frac{1}{2}c_2\pm \frac{1}{2}\sqrt{(1+3 c_2)(1-c_2)}\right.,\\ && \hfill{\left. -\frac{1}{3}\leq c_2\leq 1 \right\}}.
    \eeano where the relations between $a_4$ and $c_2$ can be obtained by considering the unit norm condition of the rows and columns.

Then observe that $HMH=\bmatrix{1 & 0\\ 0 & M_1}$ if $M\in \mathcal{C}_1$ for some matrix $M_1\in \mathcal{\overline{C}}_1,$ and  $HNH=\bmatrix{1 & 0\\ 0 & N_2}$ if $N\in \mathcal{C}_2$ for some matrix $N_2\in\mathcal{ \overline{C}}_2,$ where

  
 \beano
    \overline{\mathcal{C}}_1&=&\left\{\bmatrix{-\frac{1}{2}-a_4-c_2&-\frac{1}{2}+a_4&c_2\\c_2&-\frac{1}{2}-a_4-c_2&-\frac{1}{2}+a_4\\-\frac{1}{2}+a_4&c_2&-\frac{1}{2}-a_4-c_2}: \right. \\ && \left. a_4=-\frac{1}{2}c_2\pm \frac{1}{2}\sqrt{(1-3 c_2)(1+c_2)},-1\leq c_2\leq \frac{1}{3}\right\},\\ \overline{\mathcal{C}}_2&=&\left\{\bmatrix{\frac{1}{2}-a_4-c_2&\frac{1}{2}+a_4&c_2\\c_2&\frac{1}{2}-a_4-c_2&\frac{1}{2}+a_4\\\frac{1}{2}+a_4&c_2&\frac{1}{2}-a_4-c_2}: \right. \\ && \left. a_4=-\frac{1}{2}c_2\pm \frac{1}{2}\sqrt{(1+3 c_2)(1-c_2)},-\frac{1}{3}\leq c_2\leq 1\right\}.
    \eeano
It is to be noted that any matrix in  $\overline{\mathcal{C}}_1$  has row and column sums $-1$  while for $\overline{\mathcal{C}}_2$ it is $1.$ Also, matrices in $\mathcal{C}_1$ and $\mathcal{C}_2$ are linear combinations of at most $6$ permutation matrices.
    \end{enumerate}
Hence the proof. $\hfill{\square}$



Thus from the above theorems we obtain that: Any orthogonal matrix $A$ belonging to the spaces
$\L_i,i\in\{1,\ldots,5\}; \L_i\oplus\L_j,i,j\in\{1,\ldots,5\};\L_i\oplus\L_j\oplus\L_k,\,i,j,k\in\{1,\ldots,5\}$, $(i,j,k)\neq(1,2,5)$ and $\L_2\oplus\L_3\oplus\L_4\oplus \L_5,$ is always permutative or there exists $P,Q\in \P_4$ such that $PAQ$ or $H(PAQ)H$ is direct sum of OPMs. Finally we conclude this section with the following important remark about orthogonal matrices of order $4$ that are direct sum of OPMs. 

\begin{remark}{(Orthogonal matrices in $\L$ that are direct sum of OPMs)} Let $A\in \L$ be a  $4 \times 4$ orthogonal matrix such that $PAQ$ is the direct sum of OPMs for some $P,Q\in \P_4$. Then $PAQ = \bmatrix{1&0\\0&B}$ for $B\in \overline{\X}_1\cup \overline{\Z}_1,$ or $PAQ = \bmatrix{-1&0\\0&C}$ for $C\in \overline{\Y}_{-1}\cup \overline{\mathcal{W}}_{-1},\,$ where
 \begin{eqnarray*}
\overline{\X}_1 &=& \left\{ \bmatrix{x & y &1-x-y\\1-x-y & x & y\\y &1-x-y & x} : x^2+y^2-x-y+xy=0 \right\},\\
\overline{\Y}_{-1} &=&  \left\{ \bmatrix{
x & y &-1-x-y\\-1-x-y & x & y\\y &-1-x-y & x} : x^2+y^2+x+y+xy=0 \right\},\end{eqnarray*}
 $\overline{\Z}_{1}=\{\overline{P}A : A\in\overline{\X}_1\},\overline{\mathcal{W}}_{-1} = \{\overline{P}B : B\in \overline{\Y}_{-1}\}$ and $\overline{P}$ is the $3 \times 3$ permutation matrix corresponds to the permutation $(23).$ Note that union of $\overline{\X}_1, \overline{\Y}_{-1}, \overline{\Z}_1$ and $\overline{\mathcal{W}}_{-1}$ provides the set of all permutative orthogonal matrices OPMs of order $3$ \cite{sarkar2020}. Obviously the above matrices $PAQ$ are not permutative. Clearly if $PAQ$ is the direct sum of two $2\times 2$ permutative matrices, then $PAQ\in \P_4$ or equivalently $A\in \P_4.$



\end{remark}
 

 

\noindent{\bf Conclusion.} In this paper, we have provided parametric representation of all orthogonal permutative matrices (OPMs) of order $4$ over the field of complex, real and rational numbers. Consequently, we have shown that OPMs can be written as linear combinations of permutation matrices. We have determined several matrix spaces such that any orthogonal matrix $A$ in these spaces is always permutative  or $PAQ$ or $HPAQH$ is a direct sum of OPMs for some permutation matrices $P, Q$ and the Hadamard matrix $H$ given in equation (\ref{hadamard}). These matrix spaces are defined by direct sums of linear spaces $\L_i, i=1,\hdots,5$ which are generated by linearly independent Hadamard orthogonal permutation matrices. However, exploring all such matrix spaces and characterizing combinatorial structure of all
orthogonal matrices of order $4$ is beyond the scope of this paper. For example, consider the orthogonal matrix \beano M= \frac{1}{11}\bmatrix{10&-2&-1&4\\-2&7&-2&8\\-1&-2&10&4\\4&8&4&-5}&=&\frac{1}{11}P_{(12)}+\frac{7}{11}P_{(34)}-\frac{1}{11}P_{(13)(24)}+\frac{4}{11}P_{(14)(23)}\\
 &&+\frac{9}{11}P_{(24)}-\frac{3}{11}P_{(12)(34)}-\frac{6}{11}P_{(23)}\in \L_1\oplus \L_2 \oplus \L_5.\eeano It can be verified that $M$ does not have any of the combinatorial structure as mentioned above. We plan to explore this problem in future. \\



\noindent{\bf Acknowledgement.} Amrita Mandal thanks Council for Scientific and Industrial Research (CSIR), India for financial support in the form of a junior/senior research fellowship.


\end{document}